\documentclass[11pt]{article}
\usepackage{a4}
\usepackage{amssymb}
\usepackage{hyperref}

\let\ocirc=\odot
\let\nonu=\nonumber

\def\su{\ifmmode SU(2) \else $SU(2)$\fi}

\newcommand\I{\mathbb I}
\newcommand\C{\mathbb C}
\newcommand\Z{\mathbb Z}
\renewcommand\L{\mathcal L}
\newcommand\R{\mathbb R}
\newcommand\tr{\mathop\mathrm{Tr}}
\newcommand\alg{{\mathcal A}}
\newcommand\ch{\mathop\mathrm{Ch}}
\newtheorem{defin}{Definition}

\begin{document}

\title{Lectures on the three--dimensional non--commutative spheres.}
\author{Marc P. Bellon\thanks{On leave from:
Laboratoire de Physique Th\'eorique et Hautes Energies, Boite 126,
4 Place Jussieu,  75252 Paris Cedex 05. Unit\'e Mixte de Recherche
UMR 7589, Universit\'e Pierre et Marie Curie-Paris6; CNRS;
Universit\'e Denis Diderot-Paris7.}\\
 {\normalsize \it CEFIMAS, Av.\,Santa Fe 1145,
 1069 Capital Federal, Argentina}\\
\normalsize\it Departamento de F\'\i sica,
Universidad Nacional de La Plata\\ {\normalsize\it C.C. 67, 1900
La Plata, Argentina}}
\date{}
\maketitle
\begin{abstract}
These are expanded notes for a short course given at the Universidad
Nacional de La Plata. They aim at giving a self-contained account of
the results of Alain Connes and Michel Dubois--Violette.
\end{abstract}

\section{Introduction}

The description of the moduli space of non--commutative three--dimensional
spheres is a beautiful piece of mathematics, first presented by
Alain Connes and Michel Dubois--Violette in~\cite{CDV1,CDV2,CDV3}. In
non--commutative geometry, it is the first exhaustive study of a category
of non--commutative spaces.

Non--commutative spaces have been studied both from the mathematical and
physical perspectives. In the mathematical realm, a strong motivation
has been the study of ``bad" quotients of spaces, like foliations with
dense leaves and quotients by ergodic group actions.

Physical motivations have come from the realm of quantum gravity, since
absolute limitations on the measurements of positions of events might
translate in non-commuting coordinates.  In string theories, background
gauge field may induce non--commutation of position fields.

The most complete presentation of non-commutative geometry can be found
in the book of the creator of the field, Alain Connes~\cite{Co88},
but the shorter presentation of~\cite{Kha04} provides an interesting
introduction. Finally, a panorama of the domain and its applications
can be found in~\cite{CoMa06}.

I will first briefly motivate the definition of non-commutative
spheres, showing in particular that ordinary spheres as well as some
non-commutative deformations are included. In one and two dimensions,
the only objects which satisfy the conditions are the commutative spheres,
so that we have to go to three dimensions to have non-trivial solutions.
Then the moduli space of three-dimensional non-commutative spheres
is introduced, the corresponding quadratic algebras are described. A
large place is devoted to the geometric data, which allows to show that
differing parameters really describe differing algebras, with some lights
about the special cases.   The last sections introduce the covering of the
generic quantum spheres by a cross-product algebra of functions with $\Z$.

These lecture notes use two principles to reach a fairly complete
treatment in a reasonable space. Avoid any use of the transcendental
functions and only present one of the different possible presentations. In
particular, I always work with the variables which are used only for the
complexified moduli space in~\cite{CDV3}. This means that the generator
of the algebra are not self-adjoint and it requires some care to express
the reality conditions, but it greatly simplifies other aspects. 

\section{Basic facts in non--commutative geometry}

In ordinary differential geometry, a large use is made of the differential forms
which allow to define the De Rham complex. Forms can be considered as elements of
an algebra generated by the functions and symbols $df$ for any function $f$ with
the relations \[d1=0,\quad d(fg) = (df) g + f(dg),\quad \hbox{(Leibniz rule).}\]
The differential $d$ can be extended to these forms by making it obey a graded
Leibniz rule and imposing that $d^2=0$. The classical De
Rham complex is obtained adding the commutation rules
\[f (dg) = (dg) f,\quad df\; dg = - dg\; df\]
However these rules are inconsistent if the product of $f$ and $g$ is not
commutative. An important step has therefore been to realize that cyclic symmetry
is sufficient to establish many properties, giving rise to cyclic (co-)homology.

A fundamental object for any non--commutative space $\alg$ is the K-homology, which
takes different forms in even and odd dimensions. The even dimensional case is
linked to the classification of finitely generated projective modules. Projective
modules over the algebra of functions can be realized as the vector space of
sections of vector bundles and therefore allow to define direct analogs of
K-groups. The trivial modules corresponding to trivial bundles are simply powers
$\alg^n$. The projective modules are defined by projection operators $e$ in
$M_n(\alg)$ singling out a submodule of $\alg^n$. Instead of the projector $e$,
which satisfy $e^*=e$ and $e^2=e$, it can be easier to deal with the involution
$s=2e-1$ which is also hermitian and satisfy $s^2=1$.

In the odd dimensional case, the generators of K-theory are unitary operators $U$
in $M_n(\alg)$. They are characterized by $UU^*=U^*U=1$. In both case there is a
coupling between the K-groups and the cyclic homology, the Chern characters. To
describe them, it is convenient to introduce the operation $\ocirc$ in
$M_n(\alg)\simeq\alg\otimes M_n(\mathbb C)$:
\begin{eqnarray}
\ocirc:&& M_n(\alg)\times M_n(\alg) \rightarrow M_n(\alg\otimes\alg) \nonu\\
	&& (a\otimes m, b\otimes n) \rightarrow (a\otimes m)\ocirc (b\otimes n) =
			(a \otimes b) \otimes (m\circ n)
\end{eqnarray}
With a trace acting only on the matrix part of the expressions, we have the
following formulas for the Chern characters:
\begin{eqnarray}
\mathrm{Ch}_{2n}(s) &=& \mathrm{Tr}( \overbrace{s\ocirc s\ocirc \cdots \ocirc
s}^{2n+1} ) \\
\mathrm{Ch}_{2n+1}(U,U^*) &=& \mathrm{Tr}(
\overbrace{\widehat{U\ocirc U^*} \ocirc \cdots \ocirc \widehat{U\ocirc U^*}}^{n+1}
\nonu \\ && \qquad
-\overbrace{\widehat{U^*\ocirc U} \ocirc \cdots \ocirc \widehat{U^*\ocirc U}}^{n+1})
\end{eqnarray}

\section{Definition of non--commutative spheres.}

The sphere of dimension $n$ has the property that the cohomology and the homology is
supported only in the dimensions 0 and~$n$. Equipped with the notion of Chern
character in cyclic cohomology, the following definition has been
proposed~\cite{CDV1} for the non--commutative spheres.
\begin{defin}
An even non-commutative sphere of dimension $2n$ is an algebra generated by the
matrix elements of a hermitian involution $s$ such that $\ch_{2p}(s)=0$ for all
$p$ smaller than $n$ and $\ch_{2n}(s)$ is the volume form of the sphere.

An odd non-commutative sphere of dimension $2n+1$ is an algebra generated by the
matrix elements of a unitary $U$ together with it hermitian conjugate $U^*$
such that $\ch_{2p+1}(U,U^*)=0$ for all
$p$ smaller than $n$ and $\ch_{2n+1}(U,U^*)$ is the volume form of the sphere.
\end{defin}
In many cases, it will be convenient to separate the conditions that $U$ is
unitary or $s$ is an involution in a first one that $s^2$ or $UU^*=U^*U$ be
proportional to the identity
\begin{equation}
s^2 = C\I, UU^*=U^*U = C\I
\end{equation}
and a second one that this proportionality constant
is 1. The advantage is that the first part is purely quadratic in the matrix
elements of $s$ or $U$ and $U^*$. The proportionality constant $C$ can be seen
to be central in the
algebra by considering the products $s^3$, $U U^*U$ or $U^* U U^*$. It is
therefore always possible to consider the quantum sphere as a quotient of a
quantum affine space which is a quadratic algebra by the ideal generated by
$C-1$, where $C$ is a central element quadratic in the generators.

With this formula, it is easy to see that a quantum sphere allows to define a
quantum sphere with one dimension more, by an operation called suspension.
If we add a central hermitian
element $x$ to the algebra, $U=x\I+is$ and $U^*=x\I-is$ satisfy
$UU^*=U^*U=(x^2+C)\I$, with the new constant $C'=x^2+C$. The odd Chern forms are
cup products of the even Chern forms with the homology of $\C[x]$, so that if
$s$ defines a sphere of dimension $2n$, $U$ and $U^*$ define a sphere of
dimension $2n+1$.  If we start from an odd sphere, we have to double the matrix
dimension to define $s$ through:
\begin{equation}
s = \pmatrix { x\I & U\cr U^*& -x\I\cr}
\end{equation}
The quadratic relation is satisfied again with $C'=C+x^2$, $\tr s$ is
evidently 0 and the Chern forms of higher rank are again cup products
of the odd Chern forms for $U$ with the ones of the line.

In~\cite{CDV1}, it was shown that the commutative
spheres, as well as some non--commutative deformations, satisfy this
definition. Indeed, the matrix $s$ or $U$ can be defined in terms of
the Clifford algebra of the space of dimension $2n+1$:
\begin{eqnarray}
s &=& \sum_{j=0}^{2n} x_j \gamma_j, \\
U &=& x_0 \mathbb I + i \sum_{i=1}^{2n+1} x_j \gamma_j.
\end{eqnarray}
The properties of the Clifford algebra ensure that
\begin{eqnarray}
  s^2 &=& \sum_{j=0}^{2n} x_j^2 \mathbb I ,\\
  U U^* &=& U^* U =  \sum_{j=0}^{2n+1} x_j^2 \mathbb I,
\end{eqnarray}
so that the defining properties of the involutions or the unitary operators
are satisfied on the sphere. The properties of the Clifford algebra allow to
show that the Chern characters of low order are. In the even case,
the Chern character involves the
trace of an odd number of $\gamma$ matrices. This is zero unless the
number of matrices is at least equal to $2n+1$, in which case it is
proportional to the totally antisymmetric symbol. The low order Chern
characters are therefore zero and we get for the character $\ch_{2n}$
that it is proportional to the volume form. In the odd case, it can be
seen that the two terms annihilate themselves if there are not an odd
number of terms proportional to $x_0$. Again we remain with the trace  of
an odd number of $\gamma$ matrices, and the lowest order non-zero Chern
character is $\ch_{2n+1}$ with a single term proportional to $x_0$. Again
it can be shown that this non-zero Chern character is proportional to
the volume form.
One could also obtain the commutative spheres by iterated suspension of the
circle, the sphere of dimension 1.

In low dimensions, the only possible spheres are the commutative ones. In
dimension one, there are no conditions and we have the algebra generated
by a single unitary element $u$ which is the algebra of function on the
circle. In dimension two, the generator $s$ must be taken as a two by
two matrix and it satisfies the condition $\ch_0(s)=0$ which is simply
that it has zero trace. We can therefore write it, allowing for the
hermitian constraint:
\begin{equation}\label{2dim}
    s = \pmatrix{ x & y\cr y^* & -x\cr}.
\end{equation}
The condition $s^2=\mathbb I$ then implies that $x$, $y$ and $y^*$
mutually commute and satisfy $x^2+y y^* =1$. The first non-trivial case is
the three--dimensional one.

\section{Moduli space of non-commutative three--spheres}

For three-dimensional spheres, the generator of K--theory is a unitary
operator $U$ in $M_2(\alg)$.  In terms of the Pauli matrices $\sigma$,
$U$ and its adjoint are written:
\begin{eqnarray}	\label{fond}
	U &=& z_0 \mathbb I + i \sum_{j=1}^3 z_j \sigma_j. \\
	U^* &=& z^*_0 \mathbb I - i \sum_{j=1}^3 z^*_j \sigma_j.
\end{eqnarray}
The defining equations of the three--sphere are invariant by left and
right multiplication of the matrix $U$ by unitary matrices
$A$ and $B$. This corresponds to an arbitrary $SO(4)$ rotation mixing
the algebra elements $z_\mu$ and their multiplication by a common phase.

With these parameters,  the Chern character $\ch_1(U,U^*)$
reads:
\begin{equation}
\sum_{\mu=0}^3 z_\mu \otimes z^*_\mu - \sum_{\mu=0}^3 z^*_\mu \otimes
z_\mu . \label{cond}
\end{equation}
For this Chern character to be zero, the $z^*_\mu$ must be linear
combinations of the $z_\mu$. If the $z_\mu$ are linearly independent,
there exists a unique matrix $\Lambda$ such that
\begin{equation}\label{lamb}
	z_\mu^* = \Lambda_{\mu\nu} z_\nu
\end{equation}
The case where the $z_\mu$ are linearly dependent is not really important
since it does not allow to form spaces with non--trivial
three--volumes. We can use the
redefinition by a four dimensional rotation to put some generators to
zero and remain with independent generators. A formula similar to
eq.~(\ref{lamb}) can be written for the non-zero generators.

Inserting in eq.~(\ref{cond}) shows that $\Lambda$ is symmetric and
taking the adjoint of eq.~(\ref{lamb}) shows that it is unitary. The
real and imaginary part of $\Lambda$ are real symmetric matrices, the
unitarity of $\Lambda$ implies that they commute and they can
therefore be simultaneously diagonalized by a four dimensional
rotation:
\begin{equation}
 z^*_\mu = \lambda_\mu z_\mu , \label{lam}
\end{equation}
without summation over $\mu$.
The four eigenvalues of $\Lambda$ are unit norm complex
numbers and they can all be multiplied by a common phase by a
redefinition of $U$.  The moduli space is therefore of dimension three.

After the diagonalization, we keep  the freedom of changing the order of
the eigenvalues.  There is a natural cyclic order on the circle, but the
determination of the starting point, which can be set to one, remains. This can
be settled if we impose that the imaginary parts of $\lambda_2/\lambda_0$ and
$\lambda_3/\lambda_1$ are positive. In terms of the parameters $\phi_i$ such
that:
\begin{equation}
\lambda_0=1,\quad \lambda_i = e^{i\phi_i},
\end{equation}
this gives  two possible set of constraints:
\begin{eqnarray}
	\label{A} 0 \leq \phi_1 \leq \phi_2 \leq \phi_3 \leq \pi, \\
	\label{B} \phi_1 \leq \phi_2 \leq \pi \leq \phi_3\leq \pi+\phi_1.
\end{eqnarray}
The boundaries of this fundamental domain correspond to either the
equality of two $\lambda$ or to a relation $\lambda_i = -\lambda_j$.

It remains to show that different parameters really
correspond to distinct algebras and that an algebra can be realized with any
choice of these parameters.

\section{The quadratic algebras.} \label{quad}
The unitarity condition on $U$ can be expanded on the basis of the
matrices $\mathbb I$ and $\sigma_i$. The $\mathbb I$ term is identical
for the $UU^*$ and $U^*U$ products and gives the quadratic relation:
\begin{equation} \label{norm}
 \sum_{\mu=0}^3 \lambda_\mu z_\mu^2 = 1.
\end{equation}
This is the inhomogeneous relation in the algebra and its left hand side
is automatically central.

There remains six equations which can be deduced from the following two by a
circular permutation of the indices 1, 2 and 3.
\begin{eqnarray}\label{relation}
\lambda_0 z_1 z_0 - \lambda_1 z_0 z_1 + \lambda_3 z_2 z_3 - \lambda_2 z_3 z_2 &=& 0,
\nonu\\
-\lambda_1 z_1 z_0 + \lambda_0 z_0 z_1 + \lambda_2 z_2 z_3 - \lambda_3 z_3 z_2 &=& 0.
\end{eqnarray}

These two equations are related if we take new
generators $z'_0= -\lambda_0 z_0$, $z'_i = \lambda_i z_i$, and define new
parameters $\lambda'_\mu= \lambda^{-1}_\mu$.
This is linked to the symmetry of the roles of $U$ and $U^*$ in the
definition of the algebra. The primed parameter determine $U^*$ in the same
manner that the other ones determine $U$. This only changes the sign of the
Chern character, that is the orientation of the sphere. We have therefore an
isomorphism between the algebra of a sphere and
the algebra for the sphere with all of the parameters $\lambda$ replaced by
their inverses. Since the $\lambda$ are unit norm complex numbers, taking the
inverses or the complex conjugates are equivalent.

The sum and the difference of the equations~(\ref{relation}) express
them in terms of commutators and anticommutators:
\begin{eqnarray}\label{comm}
(\lambda_1 - \lambda_0) [z_0,z_1]_+ &=& (\lambda_2+\lambda_3) [ z_2,z_3]_-,
\nonu \\
(\lambda_1 + \lambda_0) [z_0,z_1]_- &=& (\lambda_3-\lambda_2) [ z_2,z_3]_+.
\end{eqnarray}
In this form, it is clear that when the $\lambda$'s are equal, all these
relations reduce to the mutual commutation of all the generators $z_\mu$, so
that we obtain the commutative sphere.

We can also remark that the algebra and its opposite defined by $a^\circ
b^\circ = (ba)^\circ $ are isomorphic. Indeed, the commutators change
sign and the anticommutators are unchanged in the change to the opposite
algebra and this can be compensated in eqs.~(\ref{comm}) by changing the sign of
any of the generators, since each of them appears once in each of the
equations~(\ref{comm}). Changing the sign of two generators yields an
automorphism.

In the form~(\ref{comm}), changing the scale of the generators can allow
to put one of the equations to a standard form, if we are not in one of
the special cases where some of the sums or differences of the $\lambda$
are zero. According to the one put to standard form, we can obtain a
Sklyanin algebra or a similar one given by the equations:
\begin{eqnarray}\label{skly}
  [Z_0,Z_1]_- &=& [ Z_2,Z_3]_+ ,\nonu\\
  (\lambda_1 - \lambda_0)(\lambda_2-\lambda_3) [Z_0,Z_1]_+ &=&
  (\lambda_1 + \lambda_0)(\lambda_3+\lambda_2) [ Z_2,Z_3]_-.
\end{eqnarray}
Again, the equations with circular permutations of the indices are
understood. From now on, we will adopt the convention that the indices
$k,l,m$ will stand for $1,2,3$ or any of its circular permutations. In any
expression with these three indices, a sum over $k$ or a product over $k$
will stand for the sum or the product on the three possible permutations.
If we introduce the quantities $a_k = (\lambda_k +\lambda_0)(\lambda_l+\lambda_m)$,
the second set of equations takes the form
\begin{equation}\label{skly2}
    (a_m - a_l) [Z_0,Z_k]_+ = a_k [ Z_l,Z_m]_-.
\end{equation}
This shows that the quadratic part of the algebra only depends on
two parameters, since these equations are invariant by a rescaling of
the $a_i$. This corresponds to the fact that there is always a second
central element quadratic in the generators, so that the sphere can be
deformed by modifying the inhomogeneous condition by adding this other
central element. This second central element gives a foliation of the sphere in
tori.

The explicit relation between $z_\mu$ and $Z_\mu$ is given by
\begin{equation}
Z_\mu = \rho_\mu z_\mu,
\end{equation}
with $\rho$ satisfying the following equations:
\begin{equation}
\rho_0 \rho_k (\lambda_m - \lambda_l) = \rho_m \rho_l (\lambda_0 + \lambda_k)
\end{equation}
We can choose
\begin{eqnarray}
\rho_0^2 &=& \prod (\lambda_0+\lambda_k), \quad \rho_k^2 =
(\lambda_0+\lambda_k)(\lambda_l-\lambda_k)(\lambda_k-\lambda_m) \nonu \\
\rho_0 \rho_1 \rho_2 \rho_3 &=& \prod_k \bigl[ (\lambda_0 + \lambda_k)
(\lambda_m - \lambda_l)\bigr]
\end{eqnarray}
It is interesting to note that with this choice for $\rho_\mu$ and
the normalization $\prod \lambda_\mu = 1$, the new generator $Z_\mu$
are hermitian or anti-hermitian.
\begin{equation}
 Z_\mu^* = \rho^*_\mu z_\mu^* = \frac{\lambda_\mu \rho_\mu^*}{\rho_\mu} Z_\mu
\end{equation}
The relation does not depend on the choice of square roots in the
definition of $\rho_\mu$. One has for the unit norm $\lambda$:
\begin{equation}
\frac{\lambda_\mu+\lambda_\nu}{\lambda_\mu^*+\lambda_\nu^*}
= \frac{\lambda_\mu+\lambda_\nu}{1/\lambda_\mu+1/\lambda_\nu}
= \lambda_\mu \lambda_\nu
\end{equation}
The square of the factor between $Z_\mu$ and $Z_\mu^*$ is one. The product
of these factors however is $-1$ using the constraint on the product of
the $\rho$. This result could be anticipated by considering the adjoint
of eqs.~(\ref{skly}).  In the case of eq.~(\ref{A}), one can further
show that all the new generators are hermitian except $Z_2$.

With this presentation of the algebra, the center of the algebra contains
three elements $Q_k$, whose sum is zero:
\begin{equation}\label{center}
    Q_k = (a_m-a_l) (Z_0^2 + Z_k^2) + a_k ( Z_m^2-Z_l^2)
\end{equation}
The proof uses the identities:
\begin{equation}\label{iden}
    [Z_\mu,Z_\nu^2]_- = [ [Z_\mu,Z_\nu]_-, Z_\nu]_+ = [ [Z_\mu,Z_\nu]_+, Z_\nu]_-.
\end{equation}
By choosing the appropriate one of these two possibilities for each term
of the commutator of $Z_\mu$ and $Q_k$, the identities~(\ref{skly},%
\ref{skly2}) allow to obtain $a_k$ or $a_l-a_m$ as
common factor of double (anti-)commutators which sum up to zero.

\section{Geometric data.}
In order to show that the algebras associated to distinct parameters
$a_i$ are really inequivalent, we introduce the geometric data. The
quadratic relations $R$ defining the algebra are elements of the tensor
product $V\otimes V$ of the vector space of the generators of the
algebra. They define quadratic equations in the product of projective
spaces $PV^* \times PV^*$.  The resulting space $\Gamma$ define a
correspondence $\sigma$ between the projection of $\Gamma$ on the first and the
second factor. The association of the projection $E$ of $\Gamma$ and
the correspondence $\sigma$ defined by $\Gamma$ is the geometric data.
Isomorphic algebras will give isomorphic geometric data, so that we can
show that some quadratic algebras are inequivalent.

These geometric data, plus the line bundle stemming from the canonical one on
$PV^*$, will furthermore allow to represent the algebra on a space of sections
of line bundles.

It appears that the first form of the relations~(\ref{relation}) is the most
interesting for the determination of the geometric data. If we denote by
$(y_\mu)$ and $(y'_\mu)$ coordinates in the first and second $PV^*$ factor in
a basis dual to the $z_\mu$, the equations of $\Gamma$ read:
\begin{eqnarray}
\lambda_0 y_1 y'_0 -\lambda_1 y_0 y'_1  +  \lambda_3 y_2 y'_3-\lambda_2 y_3 y'_2
&=& 0, \nonu\\
-\lambda_1 y_1 y'_0 +\lambda_0 y_0 y'_1  +  \lambda_2 y_2 y'_3-\lambda_3 y_3 y'_2
&=& 0,
\end{eqnarray}
plus the one deduced by circular permutation of the indices (1,2,3). The point
of E are such that there exist a non zero solution in $(y'_0,y'_1,y'_2,y'_3)$.
This means that the following matrix is of rank less than four:
\begin{equation} \label{matrice}
\pmatrix{
\lambda_0 y_1 &-\lambda_1 y_0 &-\lambda_2 y_3  & \lambda_3 y_2 \cr
\lambda_0 y_2 &\lambda_1 y_3 &-\lambda_2 y_0  & -\lambda_3 y_1 \cr
\lambda_0 y_3 &-\lambda_1 y_2 &\lambda_2 y_1  & -\lambda_3 y_0 \cr
-\lambda_1 y_1 &\lambda_0 y_0 &-\lambda_3 y_3  & \lambda_2 y_2 \cr
-\lambda_2 y_2 &\lambda_3 y_3 &\lambda_0 y_0  & -\lambda_1 y_1 \cr
-\lambda_3 y_3 &-\lambda_2 y_2 &\lambda_1 y_1  & \lambda_0 y_0 \cr
}
\end{equation}
Each of the fifteen four by four determinants that can be extracted from
this matrix must be simultaneously zero. This calculation is made
manageable by the factorization of the four three by three minors extracted
from the first three lines of the matrix: they all have the common quadratic factor
$y_0^2+y_1^2+y^2_2+y_3^2$. Similarly, the three by three minors extracted from
the last three lines have the common factor $\sum \lambda_\mu^2 y_\mu^2 $.
Therefore on the curve defined by the intersection of these two quadrics, six
of the four by four submatrices have determinant 0. Furthermore, since the upper
and the lower subblock are of rank 2, it will not be necessary to verify that
the remaining nine submatrices are of determinant 0. Indeed, apart from isolated
points on the curve, the third line of each group is a linear combination of the
two others. If the determinant of the matrix formed with the lines 1,2,4,5 is
0, we therefore have that they generate a space of dimension at most 3 and all
other determinants are also zero. Continuity allows to ensure that these other
determinants remain zero on the points where the lines 1 and 2 for example are
proportional. That the determinant of the matrix formed with the lines 1,2,4,5
is zero can be seen by writing it:
\begin{equation}
	(\lambda_0 \lambda_3 + \lambda_1 \lambda_2) \left[ ( y_1^2 + y_2^2)
	( \lambda_0^2 y_0^2 + \lambda_3^2 y_3^2 )-( y_0^2 + y_3^2)( \lambda_1^2
	y_1^2 + \lambda_2^2 y_2^2 )\right]
\end{equation}
The two quadratic equations define a curve which is generically an elliptic one:
\begin{eqnarray}
 	\label{generic1}
	y_0^2 +y_1^2 +y_2^2 +y_3^2 &=&0, \\
	\lambda_0^2y_0^2 +\lambda_1^2y_1^2 +\lambda_2^2y_2^2 +\lambda_3^2y_3^2 &=&0.
	\label{generic2}
\end{eqnarray}
If different $\lambda^2$ become equal, this curve decompose in a product of
rational curves.

Suppose that equation~(\ref{generic2}) is not verified. We then have
three necessary conditions stemming from the determinants including
the lines 4,5,6.
\begin{eqnarray}\label{special}
  (\lambda_0^2-\lambda_1^2) y_0 y_1 &=& (\lambda_2^2-\lambda_3^2) y_2 y_3
  \nonu \\
  (\lambda_0^2-\lambda_2^2) y_0 y_2 &=& (\lambda_3^2-\lambda_1^2) y_3 y_1
  \nonu\\
  (\lambda_0^2-\lambda_3^2) y_0 y_3 &=& (\lambda_1^2-\lambda_2^2) y_1 y_2
\end{eqnarray}
Again, these equations only depend on the squares of the $\lambda_\mu$. In
the case where all the squares are distinct, these homogeneous equations
have exactly eight projective solutions. The four trivial solutions
$P_\mu$ such that only $y_\mu$ is non zero and four other ones which
can be deduced from one of them by
changing the sign of a pair of $y_\mu$. However these additional solutions
belong to the generic variety and we do not have to consider them independently.
Indeed, in the case where all the $y_\mu$ are not zero, the product of the first
two equations for example gives the ratio $y_0^2/y_3^2$ so that we have:
\begin{eqnarray}
 y_0^2 &=& (\lambda_2^2-\lambda_3^2) (\lambda_3^2-\lambda_1^2)
 (\lambda_1^2-\lambda_2^2) \nonu \\
 y_1^2 &=& (\lambda_2^2-\lambda_3^2) (\lambda_0^2-\lambda_2^2)
 (\lambda_0^2-\lambda_3^2) \nonu \\
 y_2^2 &=& (\lambda_0^2-\lambda_1^2) (\lambda_3^2-\lambda_1^2)
 (\lambda_0^2-\lambda_3^2) \nonu \\
 y_3^2 &=& (\lambda_0^2-\lambda_1^2) (\lambda_0^2-\lambda_2^2)
 (\lambda_1^2-\lambda_2^2)
\end{eqnarray}
The expressions $\sum y_\mu^2$ and $\sum \lambda_\mu^2 y_\mu^2$ are polynomials of
degree two in say $\lambda_0^2$ which vanish for the three values $\lambda_i^2$,
therefore they are zero. The only particular cases to consider are therefore the
four projective points with only one non-zero coordinate. It is easy to verify
that each of these points is indeed in the variety $E$ and that they are
associated to themselves by the correspondence $\sigma$.

The case where the equation~(\ref{generic1}) is not verified is completely
similar.  In fact, a change of variable allows to exchange the upper and
lower block of the matrix, so that no new calculations are necessary to
ensure that in the generic case, the variety $E$ is the union of the
elliptic curve defined by the equations~(\ref{generic1},\ref{generic2}) and the
four points $P_\mu$. The correspondence $\sigma$ is the identity on the
points $P_\mu$. The geometric data for the opposite algebra involve the
same variety $E$ but the inverse correspondence $\sigma^{-1}$. The
isomorphism of the algebras of the three-spheres with their opposites
implies that $\sigma^{-1}$ can be obtained from $\sigma$ simply by
conjugating with the change of sign of $y_0$. This identifies $\sigma$ as
a translation of the elliptic curve. We therefore have two parameters to
describe the geometric data, the modulus of the elliptic curve and the
parameter of the translation.

When different $\lambda_\mu^2$ become equal, the situation is very
similar. The equations~(\ref{generic1},\ref{generic2}) define a curve
which now split in different rational components. As in the generic case,
if one of these equations is not satisfied, the equations~(\ref{special})
or a similar set are necessarily satisfied.  The solutions  are either
included in the generic variety or have some zero components. Substituting
in the matrix allows to show that these points with some zero components
really are in the characteristic variety and that the correspondence is
either the identity or a symmetry with some of the coordinates changing
sign. The full description of these non-generic cases is better done
after the study of symmetries which allow to relate some cases.

\section{Reality and symmetries.}
Before completing the analysis of the geometric data, we have to consider
the reality conditions and the symmetries for the algebras and the relations
between different algebras that they
allow.

The generators $z_\mu$ allow for a simpler description of the
characteristic variety, but their conjugation properties make the reality
conditions for this variety subtler. The restriction to the space of
generators $V$ of the adjunction in the algebra maps to $V$, so that
it defines an antilinear map $j$. In $V^*$, the complex conjugation is
defined by $j(L)(Z) = \overline{L(j(Z))}$, which gives in coordinates:
\begin{equation}\label{realV}
(y_0,y_1,y_2,y_3)\longrightarrow( \lambda_0^*
y_0^*,\lambda_1^* y_1^*,\lambda_2^* y_2^*,\lambda_3^* y_3^*)
\end{equation}
Since the $\lambda$ are complex numbers of unit norm,
they do not obey the reality condition $x^*=x$ but the one
$x^*=x^{-1}$. Combined with the complex conjugation given
by eq.~(\ref{realV}), this allows to show that  the generic
variety is indeed invariant by complex conjugation, since the
equations~(\ref{generic1},\ref{generic2}) on the complex
conjugates are the complex conjugates of these
same equations for the initial point.

From the formula $(ab)^*=b^* a^*$, we obtain that the set of relation $R$
is invariant by the action of $j \otimes j$ followed by the exchange of
the factors. This shows that the characteristic variety is invariant by
the action of $j$ and that we have:
\begin{equation} \label{jsig}
j\circ \sigma = \sigma^{-1} \circ j
\end{equation}

The algebra automorphisms introduced in section~\ref{quad} allow to define
twisting of the non-commutative spheres which are again non-commutative
spheres with related parameters. These automorphisms come from the
diagonal transformations in $SO(4)$ and act on the
generators by changing the sign of two of them. Now, with any algebra
*-automorphism $\tau$ one can build a cross-product algebra
$\alg \times_\tau \Z$ by adjoining
to the algebra $\alg$ a unitary symbol $W$ satisfying the relation:
\begin{equation}
 	W a = \tau(a) W
\end{equation}
In the cases where the automorphism act by linear transformation on
the generators $x_i$, we can recover an algebra with the same set of
generators by considering the subalgebra of the cross-product with
generators $x_i W$. The geometric data for this new algebra can easily
be deduced from those of the original algebra. The relations  for the
new algebra are deduced from the old ones by the action of $\I \otimes
\tau^{-1}$ so that the space $\Gamma$ is changed by the action of the
identity on the first factor and the transpose of $\tau^{-1}$ on the
second. The variety $E$ is unchanged and the correspondence becomes
${}^t\tau^{-1}\circ \sigma$.

In the case of the non-commutative three sphere, $\tau^2$ is the identity
so that we can add the condition $W^2=1$. We change the matrix generators
$U$ to the generator $U W$. The quadratic relations are deduced from the
relation $U U^* = U^* U = \I$ and keep the same form. What will change
is the relation between $U$ and $U^*$.
\begin{equation}
(U W)^* = W^* U^* = W U^*= \tau(U) W.
\end{equation}
The eigenvalues $\lambda_\mu$ corresponding to the generators on which
$\tau$ has a non-trivial action will change sign.

In particular, changing the sign of $\lambda_0$ and
$\lambda_3$ and exchanging them allows to relate 
the algebras corresponding to the case of~(\ref{B}) with those described
by~(\ref{A}).

\section{Special cases.}
In the study of the geometric data, these symmetries allow to reduce the number
of different case to study. The non-generic cases can be first classified by the
number of different values of the $\lambda^2$, which can go from 1 to 4. In the
cases of two different values of $\lambda^2$, the eigenvalues can be split in
two groups of two, or one group of one and one of three.

In the case with three different values of $\lambda^2$, we can have the two
relations $\lambda_1= \lambda_2$ or $\lambda_1=-\lambda_2$ but the two can be
related by a symmetry. The generic variety stemming from
equations~(\ref{generic1},\ref{generic2}) split
in two conics. To the four special points we must add the line through the
points $P_1$ and $P_2$.

In the case with two different values of $\lambda^2$ and unequal groups, again
the different possibilities for the sign relation between the $\lambda$ can be
related by the symmetry. The special variety now contain the plane through the
three points $P_1$, $P_2$ and $P_3$. The generic variety is now embedded in this
plane, where $\sigma$ is the identity in the simplest case where three
$\lambda$ are equal.

The case of say $\lambda_1=\lambda_2$ and $\lambda_0=\lambda_3$ divide in three
cases. The two cases where either the $\lambda$ are pairwise equal or opposite
can be related by symmetry, but the third possibility is distinct. With
$\lambda_0=-\lambda_3$ and $\lambda_1=\lambda_2$, two of the equations becomes
equal, so that the algebra grows exponentially. On the geometric data, this
translates in the apparition of coarse correspondences, where some points map
to whole lines.

Finally, when all the $\lambda^2$ are equal and can therefore be taken to be 1,
there are three cases according to the number of $-1$. With all 1, we have the
commutative case, the entire projective space is the characteristic variety
and the correspondence is the identity. The case with two $-1$ is related to
the first one by a symmetry. The case with one $-1$ is very special. All the
$a_i$ are zero and we only have three relations in the algebra. The
matrix~(\ref{matrice}) has only three independent lines and is therefore of
maximum rank 3. However, on the surface with equation~(\ref{generic1}), its
rank becomes 2 and the correspondence $\sigma$ associates with such points
whole lines.

\section{The fundamental elliptic curve.} \label{elliptic}
In section~\ref{quad}, we saw that the quadratic algebra $\R^4_\lambda$
only depends on two parameters, with all parameters lying on a curve
describing a unique quadratic algebra. The identification is linked
to a rescaling of the generators such that the algebra is defined by
eqs.~(\ref{skly},\ref{skly2}), which only depend on three homogeneous
parameters. What is most remarkable is that the curve in the space of
parameters $\lambda$ which maps to a given set of the $a$ parameters is
identical with the generic part of the characteristic variety.  This was
first an observation on the identity of the elliptic parameters, but
the two curves can be explicitly identified.

The scaled generators of the algebra are linked to a new
coordinate system for $V^*$ such that $Y_\mu = \rho_\mu y_\mu$. The
equations~(\ref{generic1},\ref{generic2}) become, after eliminating
the denominators:
\begin{eqnarray}\label{genscale}
  \prod_k (\lambda_m - \lambda_l) Y_0^2 +
		  \sum_k (\lambda_0+\lambda_m) (\lambda_0+\lambda_l)
		  (\lambda_m - \lambda_l) Y_k^2 &=& 0 , \nonu\\
  \lambda_0^2 \prod_k (\lambda_m - \lambda_l) Y_0^2 +
		  \sum_k \lambda_k^2 (\lambda_0+\lambda_m) (\lambda_0+\lambda_l)
		  (\lambda_m - \lambda_l) Y_k^2 &=& 0.
\end{eqnarray}
The important property of these equations is that they are linear
combinations of $B_k = Y_0^2 - Y_k^2$. Indeed the sums
of the coefficients are polynomials of degree less than 2 in $\lambda_0$
which are zero at the three points $-\lambda_k$.  The linear combination of
the equations~(\ref{genscale}) such that the coefficient of say $Y_3^2$
is zero allows to obtain the ratio of $B_1$ and $B_2$ on the curve. It
turns out to be the same as the one of $a_1$ and $a_2$.

The $B_k$ can be converted to the $a_k$ through a linear change of
variables with an involutive matrix $M$:
\begin{equation}\label{M}
    M = \frac{1}{2} \pmatrix{ 1 & 1 & 1 & 1\cr 1 & 1 & -1 & -1\cr
        1 & -1 & -1 & 1\cr 1 & -1 & -1 & 1\cr}.
\end{equation}
We therefore have identified the two elliptic curves, in the parameter
space of the non-commutative spheres and in the characteristic variety.

Another interest of this presentation of the quadratic algebra is that the
correspondence $\sigma$ has a simple form, since the first three relations of
the algebra, which are independent on the parameters, are sufficient to
determine the image of a point. The image $\sigma(Z)$ of a point is
determined by the condition that its image by the matrix $N(Z)$ is zero:
\begin{equation}
 N(Z) = \pmatrix{ Y_1& -Y_0& Y_3 & Y_2 \cr
				Y_2 & Y_3 & -Y_0 & Y_1 \cr
 				Y_3 & Y_2 & Y_1 & -Y_0 \cr
				(a_2-a_3) Y_1& (a_2-a_3) Y_0& -a_1 Y_3 & a_1 Y_2 \cr
				(a_3-a_1) Y_2 & a_2 Y_3 & (a_3-a_1)Y_0 & -a_2 Y_1 \cr
 				(a_1-a_2) Y_3 & -a_3 Y_2 & a_3 Y_1 & (a_1-a_2)Y_0 \cr }.
\end{equation}
The first three lines form a system of maximal rank, so that they are
sufficient to determine $\sigma(Z)$. Furthermore, if we replace $Y_0$ by
its opposite, the first three lines of $N(Z) M$ read:
\begin{equation}
	\pmatrix{ \lambda_0 &\lambda_1 &-\lambda_2 &-\lambda_3 \cr
\lambda_0 &-\lambda_1 &\lambda_2 &-\lambda_3 \cr
\lambda_0 &-\lambda_1 &-\lambda_2 &\lambda_3 \cr }.
\end{equation}
This shows that $\sigma$ is the product $I\circ I_0$ of two involutions:
$I_0$ changes the sign of $Y_0$ and $I$ corresponds to taking the inverse
of all coordinates in the $\lambda$ coordinates. These two operations are
symmetries for the elliptic curves since they have four fixed points.

On these curves there are three different notions of reality. There is the
simple one corresponding to taking real coordinates. The one associated to
the $\lambda$, which are of unit norm, so that we have $\bar v = I(v)$. As
$I$ corresponds to taking the opposite in the elliptic curve, they are
imaginary points.
Finally, the $j$ operation on $V^*$, linked to the adjunction properties
of the algebra, changes the sign of an odd number of the complex
conjugated coordinates: the $j$-real points are again purely imaginary points
of the curve, but starting from a different origin.  The transformation
$\sigma$ has real, even integer, coefficients, so that it is a translation
by a real element of the curve.

\section{Uses of the geometric data.}

The geometric data allow for the definition of a representation of the
algebra. Indeed~\cite{ATvdB}, the algebra  maps to a cross-product algebra
of sections of line bundles by the map $\sigma$. Since the variety $E$
is given as a subvariety of the projective space $PV^*$, it is equipped
with  a naturally defined line bundle $\mathcal L$, the restriction of the
canonical bundle of $PV^*$. To each element of $V$, we can associate a
section $s_V$ of this bundle $\mathcal L$. The generator of the algebra
are mapped to products $s_V W$, with the symbol $W$ satisfying the
commutation relation $W s_V = s_V \circ \sigma W$. The product of $n$
generators is therefore the product of a holomorphic section of a bundle
$\mathcal L_n$ by $W^n$. The bundle $\mathcal L_n$
is defined in terms of pullbacks of $\mathcal L$ by the maps $\sigma^n$:
\begin{equation}\label{bundle}
    \L_n := \L \otimes \L^\sigma \otimes \cdots \otimes \L^{\sigma^{n-1}}
\end{equation}
The defining relations of the algebra are satisfied from the definition
of $\sigma$.  However, this cannot provide a faithful representation for
the non-commutative spheres, since all central elements of the algebra
are represented by 0 in the generic case. Indeed, if a central element
would have a non-zero representation, we would get inconsistencies
for the divisor of its product with any generator. The trouble can
be spotted also from the dimensions of the spaces of products of $n$
factors. The Riemann-Roch theorem shows that, in the generic case, this
dimension is at most $4n$  when the dimension of the corresponding space
is $\pmatrix{n+3\cr 3\cr}$ for the quadratic algebra and $(n+1)^2$ for
the non-commutative sphere.  An other limitation of this representation
is that it is not a $*$-representation.

Alain Connes and Michel Dubois-Violette introduced an extension of this
construction based on the notion of ``central quadratic form''. This
allows for a representation where the inhomogeneous equation of the
algebra~(\ref{norm}) is satisfied. The main point is to turn the algebra
of the previous representation in a *-algebra. Taking the complex
conjugate of holomorphic functions gives antiholomorphic ones.  We
therefore introduce a second variable to represent the complex
conjugate. In a first time,  the two variables are better kept
independent, so that the space we consider is the product $E\times E$
with the complex conjugation $\tilde j$:
\begin{equation}
  \tilde j ( Z, Z') = (j(Z'),j(Z)). \label{jt} 
\end{equation}
In the end, we will restrict to the invariant subspace of $\tilde j$,
formed of points $(Z,j(Z))$.
With $c$ denoting complex conjugation in the image space, the
adjunction for the function $F$ is defined by:
\begin{equation}
	F^* = c \circ F \circ \tilde j.
\end{equation}
With the following representation for the generators,
\begin{equation} \label{repres}
\rho(Y) = Y(Z) W + W^* Y(Z'),
\end{equation}
it is easy to see that 
\begin{equation} \label{herm}
	\rho(Y)^*= \rho( j(Y) ).
\end{equation}
The transformation $\sigma$ must be extended to $E \times E$. In view
of the relation between $\sigma$ and $j$~(\ref{jsig}), we define
$\tilde \sigma$ so that it commutes with $\tilde j$:
\begin{equation}
  \tilde \sigma (Z,Z') = (\sigma(Z),\sigma^{-1}(Z'). \label{sigt}
\end{equation}
It is natural to suppose that $W$ is a unitary, so that the commutation
with $W^*$ gives a translation by $\tilde\sigma^{-1}$.

Finally, we want the map~(\ref{repres}) to define a representation of the
algebra of the non-commutative sphere. Terms $W^* Y_i(Z') Y_j(Z) W$ and
$Y_i(Z) W W^* Y_j(Z')$ will correspond to fibers in differing points of
$E$, so that a trivialization is necessary for their comparison. Such a
trivialization
will also be necessary to allow for the inhomogeneous equation~(\ref{norm}).

The product of the line bundle $\L$ on the variable $Z$ and $Z'$ is
trivialized by dividing by a quadratic form $Q(Z,Z')$. Then, the
quadratic relations read:
\begin{eqnarray}
\omega_{ij} \rho(Y_i) \rho(Y_j) &=& \omega(Z,\sigma(Z))W^2 +
(W^*)^2\omega(\sigma^{-1}(Z'),Z') \nonu \\
&&+ \frac{\omega(Z,Z')}{Q(Z,Z')}
+\frac{\omega(\sigma(Z'),\sigma^{-1}(Z))}{Q(\sigma^{-1}(Z),\sigma(Z'))}
\label{omega}
\end{eqnarray}
The parts proportional to $W^2$ or $(W^*)^2$ are zero from the defining
property of $\sigma$. The remaining part is zero
if $Q$ is symmetric and satisfies the defining property of a central
quadratic form:
\begin{equation}\label{central}
    \omega(Z,Z') Q(\sigma(Z'),\sigma^{-1}(Z)) +
    \omega(\sigma(Z'),\sigma^{-1}(Z))Q(Z,Z') =0
\end{equation}
In the generic case, this condition can be verified if both $Z$ and $Z'$
belong to the generic curve, for any quadratic form deduced from a central
quadratic elements of the algebra. The calculation is done in the
coordinates introduced in section \ref{elliptic}, for any of the central forms
$Q_k$: the transformation $\sigma$ has a simple explicit expression and the
fact that $Z$ and $Z'$ are on the curve can be expressed simply by
substituting $y_0^2 + C a_k$ to $y_k^2$. The equation~(\ref{central}) is
not verified if $Z$ is a special point and $Z'$ is on the elliptic curve,
but this does not preclude the possibility to obtain a representation of
the algebra. By specializing eq.~(\ref{central}) to the case
$Z'=\sigma(Z)$, we can further prove that $Q(Z,\sigma(Z))$ is zero is
$\sigma^4$ is not the identity. The calculation of~(\ref{omega}) applied
to $\omega = \frac12 Q$ now shows that $\frac12 Q_{ij} Y_i Y_j$ is
represented by 1. This representation is a *-representation if $Q$
verifies $Q(\tilde j(Z,Z')) = \overline{Q(Z,Z')}$.

If we now specialize to a $\tilde j$ invariant subset, i.e., to points of
the form $(Z,j(Z))$, we can further demand that $Q$ is strictly positive,
so that it defines an hermitian metric on $\mathcal L$.
This is certainly the case for the central element $C$ of the algebra,
which can also be written $\sum z_\mu z^*_\mu$. 
Division by $Q(Z,j(Z))$
does not introduce any singularity and the algebra of the sphere is
faithfully represented in the space of regular sections.
We have therefore obtained a representation of the non-commutative spheres
for generic parameters. The fact that an elliptic curve is topologically
a torus can be used to further understand the structure of the algebra.
Essentially, we have a family of non-commutative torus with a parameter
determined from the transformation $\sigma$.

However, this does not
provide for a full representation in every cases. In the case where three
of the $\lambda$ are equal, we have a central element $z_0$ which play the
role of the Planck constant for the $su(2)$ algebra formed by the three
others generators. The correspondence $\sigma$ is the identity, so that
the algebra obtained from the geometric data is commutative. This is
compatible with the algebra structure because the
characteristic variety is the union of a point and the plane $y_0=0$: on
the plane, the geometric data represent $z_0$ by 0 and all other
generators are zero on the isolated point.

The full structure of this algebra is much more interesting. We get a discrete
structure, because for a $n$ dimensional representation of the $su(2)$
algebra, $z_1^2+z_2^2+z_3^2 = (n^2-1)z_0^2$. We therefore have a discrete
family of fuzzy spheres which converge to the ordinary sphere obtained for
$z_0=0$.

\section{Conclusion}
We arrive at the end of this tour of non-commutative three-spheres.
The definition of the non-commutative three-spheres has been introduced,
the moduli space discovered.  A large place has been given to the geometric
data, first to provide for intrinsic parameters to single out the quadratic
algebras, then to show how they give rise to representations for the
algebra of the non-commutative spheres.

In a sequel to these lectures, it would be interesting to evaluate the
volume form provided by the Chern form of order $3$ and the Jacobian of
the transformation to the cross-product algebra $F_u \times_{\sigma,\L} \Z$.
In keeping to the spirit of these lectures, we could make a purely
algebraic version of the computation of section~12 of~\cite{CDV3}.

We also would like to build spectral triples for these algebras, in order
to complete the description of their geometry. However, one expects that
this would be rather difficult. In the cases where the characteristic
variety is rational, which have quantum group structures of $SU_q(2)$, the
construction has met with some difficulties: in the generic case, we do
not have any symmetry to assist in the process.


\begin{thebibliography}{1}

\bibitem{CDV1}
A.~{Connes} and M.~{Dubois-Violette}.
\newblock {Noncommutative Finite-Dimensional Manifolds. I. Spherical Manifolds
  and Related Examples}.
\newblock {\em Communications in Mathematical Physics}, 230:539--579, 2002.

\bibitem{CDV2}
Alain Connes and Michel Dubois-Violette.
\newblock Moduli space and structure of noncommutative 3-spheres.
\newblock {\em Lett. Math. Phys.}, 66:91--121, 2003.
\newblock arXiv:math/0308275.

\bibitem{CDV3}
A.~{Connes} and M.~{Dubois-Violette}.
\newblock Noncommutative finite-dimensional manifolds. {II}. {M}oduli space and
  structure of non commutative 3-spheres.
\newblock 2005.
\newblock arXiv:math/0511337.

\bibitem{Co88}
Alain Connes.
\newblock {\em Noncommutative geometry}.
\newblock Academic Press, San Diego, 1994.
\newblock available online through http://www.alainconnes.org/downloads.html.

\bibitem{Kha04}
Masoud Khalkhali.
\newblock Very basic noncommutative geometry. \hfil\break
\newblock arXiv:math/0408416, 2004.

\bibitem{CoMa06}
Alain Connes and Matilde Marcolli.
\newblock A walk in the noncommutative garden.
\newblock arXiv:math/0601054, 2006.

\bibitem{ATvdB}
M.~Artin, J.~Tate, and M.~Van den Bergh.
\newblock Some algebras associated to automorphisms of elliptic curves. the
  {G}rothendieck festschrift. vol. {I}.
\newblock {\em Prog. Math.}, 86:33--85, 1990.

\end{thebibliography}
\end{document}